# Multi-objective LQR with Optimum Weight Selection to Design FOPID Controllers for Delayed Fractional Order Processes

Saptarshi Das[a,b *], Indranil Pan[b,c], and Shantanu Das[d]

a) School of Electronics and Computer Science, University of Southampton, Southampton SO17 1BJ, United Kingdom.

b) Department of Power Engineering, Jadavpur University, Salt Lake Campus, LB-8, Sector 3, Kolkata-700098, India.

c) Department of Earth Science and Engineering, Imperial College London, Exhibition Road, London SW7 2AZ, United Kingdom.

d) Reactor Control Division, Bhabha Atomic Research Centre, Mumbai-400085, India.

**Authors' Emails:**

s.das@soton.ac.uk, saptarshi@pe.jusl.ac.in (S. Das*)

i.pan11@imperial.ac.uk, indranil.jj@student.iitd.ac.in (I. Pan)

shantanu@barc.gov.in (Sh. Das)

**Corresponding author's phone number:** +44(0)7448572598

**Abstract:**

An optimal trade-off design for fractional order (FO)-PID controller is proposed in this paper with a Linear Quadratic Regulator (LQR) based technique using two conflicting time domain control objectives. The deviation of the state trajectories and control signal are automatically enforced by the LQR. A class of delayed FO systems with single non-integer order element has been controlled here which exhibit both sluggish and oscillatory open loop responses. The FO time delay processes are controlled within a multi-objective optimization (MOO) formulation of LQR based FOPID design. The time delays in the FO models are handled by two analytical formulations of designing optimal quadratic regulator for delayed systems. A comparison is made between the two approaches of LQR design for the stabilization of time-delay systems in the context of FOPID controller tuning. The MOO control design methodology yields the Pareto optimal trade-off solutions between the tracking performance for unit set-point change and total variation (TV) of the control signal. Numerical simulations are provided to compare the merits of the two delay handling techniques in the multi-objective LQR-FOPID design, while also showing the capability of load disturbance suppression using these optimal controllers. Tuning rules are then formed for the optimal LQR-FOPID controller knobs, using the median of the non-dominated Pareto solution to handle delays FO processes.

## 1. Introduction

Modern optimal control theory has a rich set of analytical tools to design control strategies satisfying desirable characteristics of the excursion of the system states according to the designer's specifications in an optimal manner [1]. The LQR is one such design methodology whereby quadratic performance indices involving the control signal and the state variables are minimized in an optimal fashion. Historically, in the area of industrial process control, PID controllers are tuned by minimizing a suitably chosen performance index for the control loop error function, which has yielded several thousands of tuning rules [2], in order to get an optimal PID setting. The tuning of PID controller uses the knowledge of the process model (mostly integer order models) like the gain ($K$), time-constant or lag ($T$) and time-delay ($L$). In spite of the huge advancements in the theoretical aspects of optimal control, successful integration of modern optimal control techniques in practical PID control problems was not there for decades due to several hidden heuristics in the design. For example, an effective choice of the weighting matrices ($Q$ and $R$) in the optimal state feedback (LQR) design which is often impossible to know a priori, especially for the control of large industrial processes [3]. There have been some previous efforts to merge the PID controller tuning problem with LQR theory as described in [4], [5], considering the error and integral of error as the state variables. The LQR technique has also been extended for tuning PID controllers for sluggish over-damped second order processes in [5] by cancelling one of the real system poles with one of the zeros of PID controller. Thus, the approach, presented in [5] does not give the flexibility of tuning oscillatory processes by selecting the optimal controller gains *via* LQR for the three state variables i.e. error, its rate and integral. In the present paper, this concept is extended by simultaneously considering all the three state as the proportional, integral and derivative action of the controller and finding a synergism of the fractional calculus based enhancements of PID controllers [6] to circumvent the afore-said problems. The goal of the paper is to find out an answer to the following research questions – 1) how to optimally choose LQR weights, keeping in mind the final closed loop performance of a sluggish/oscillatory higher order process in FO template, 2) how to handle the time delay terms, especially large delays in LQR formulation, while preserving both the stability and performance, 3) which time delay handling technique yields a better trade-off control design in terms of Pareto non-dominance for oscillatory and sluggish higher order processes with varying level of lag and delay.

The rest of the paper is organized as follows. In section 2 we briefly introduce the background of FO multi-objective control and FO optimal control, thereby highlighting the motivation of the current research. Section 3 discusses about the optimal state feedback approach of $PI^\lambda D^\mu$ controller tuning and the methodology for the selection of LQR weighting matrices to handle FO systems with time delay. Section 4 presents the simulation studies of the proposed controller with an oscillatory and a sluggish process having small time delay, followed by test of its robustness and then verifying the methodology on lag/delay dominant and balanced lag-delay FO processes. The paper ends with the conclusions as section 5, followed by the references.

## 2. Background and motivation

### *2.1. FOPID controller and NIOPTD process*

Fractional calculus, although being an age old topic in mathematics, has only recently flourished within the domain of systems and control theory [7], [8]. The FOPID or $PI^\lambda D^\mu$ controller, proposed by Podlubny [6] is an extension of conventional PID controller and is gradually getting importance in various industrial process control applications. Due to its higher degrees of freedom, the fractional order $PI^\lambda D^\mu$ controller has better ability to enforce several conflicting control objectives than the conventional integer order PID controllers. However, the performance of such controller greatly depends on its tuning methodology [9]. Several tuning philosophies have been proposed to design PI$^\lambda$D$^\mu$ controllers e.g. analytical rule based [10], [11], stabilization based [12], time domain [9] and frequency domain [13], [12] methods, to name a few amongst other available techniques [14]. In this paper, the LQR formulation is used for tuning FOPID controllers to handle a class of FO processes with one non-integer order element ($\alpha$) as also studied in [12][15], with an additional inclusion of a time delay term along with the fractional dynamics. The reason behind the consideration of this particular FO template is that it has been shown in [9], [14], [16–18] that many higher order systems can be reduced to the Non-Integer Order Plus Time Delay (NIOPTD) template which is capable of faithfully capturing the oscillatory or sluggish higher order process dynamics with only four process parameters ($K$, $L$, $T$ and $\alpha$), thus enabling a compact representation of higher order processes.

*2.2. Determination of LQR weights in PID/FOPID controller design*

It is also well-known that an optimal state feedback regulator (LQR) automatically minimizes the variation in the state trajectories but it does not always show acceptable closed loop response and might often include high overshoot, oscillations etc. for a bad choice of the weighting matrices. In order to achieve efficient tracking for a set-point change, the weighting matrices should be chosen in such a manner that it meets some additional time domain optimality criteria in terms of overshoot, rise and settling time etc. He *et al.* [5] proposed the technique to find out the LQR weights from closed loop damping and frequency specifications. In this paper, this concept is extended with an MOO based approach to find out the optimum set of weighting matrices for the optimal regulator design, as also studied in [19], [20], [21] for the standard PID controller. Poodeh *et al.* [22] used genetic algorithm to find weighting matrices by the minimization of a custom cost function of steady-state error, maximum percentage of overshoot, rise time and settling time. Here, these concepts are extended with an LQR based framework for fractional systems with delay where the LQR weights ($Q$ and $R$) that determines FOPID gains ($K_p$, $K_i$, $K_d$) and the integro-differential orders ($\lambda$, $\mu$) of the FOPID controller both have been taken as the decision variables of the Non-dominated Sorting Genetic algorithm-II (NSGA-II) algorithm. The multi-objective design makes an optimal choice of the LQR parameters and the fractional order elements of the controller using an optimal design trade-off between two time domain performance indices, i.e. the Integral of Time Multiplied Squared Error (ITSE) and total variation of control signal which is measured as the Integral of Squared Deviation of Control Output (ISDCO) [15].

*2.3. Optimal control to LQR in the context of fractional order systems*

Similar to the conventional integer order scenario, optimal control theory has been extended for FO systems by Agrawal [23], for the Euler-Lagrange equation and boundary value problems (BVPs) with FO ordinary differential equations (ODEs). Shafieezadeh *et al.* [24] have investigated the effect of adding fractional derivatives of the state variables along with the conventional optimal state feedback law using LQR. Tricaud and Chen [25], Agrawal [26], Biswas and Sen [27] formulated the fractional optimal control problem with a finite horizon quadratic performance index involving the states and control action. The formulation has been extended by Biswas and Sen [28] for free final time optimal control problems. Tangpong and Agrawal [29], Biswas and Sen [30], and Ding *et al.* [31] also proposed similar finite horizon performance index for fractional optimal control problems and derived the optimality condition for the Euler-Lagrange equations. There have been other extensions as well e.g. optimal control theory for FO discrete time systems using the Grunwald-Letnikov approach [32], [33], optimal control of distributed systems [34], [35] and derivation of the necessary condition for optimality [36]. In spite of a large number of works on Euler-Lagrange equation, variational calculus and two point boundary value problems, there have been very few works on solving infinite horizon LQR problems and derivation of Riccati equation for FO systems. The fractional LQR was first proposed by Li and Chen [37], using the Riccati differential equation, as compared to an algebraic Riccati equation (ARE) for integer order systems. In Sierociuk and Vinagre [38] for infinite horizon LQR problems, the standard Riccati equation like solutions are obtained under some assumptions, which has also been adopted in this paper.

There are also few attempts of using the concepts of optimal control, especially LQR theory, to tune FOPID controllers. For example, Saha *et al.* [39] studied LQR equivalence of dominant pole placement problem with FOPID controllers and then used a conformal mapping based approach to approximate FOPID zeros in the primary Riemann sheet with that of a PID controller. Das *et al.* [19] studied a single objective optimization based optimum weight selection of discrete time LQR to tune digital PID controllers using FO integral performance index. Das *et al.* [40] extended this concept to design an LQR-FOPID to handle FO systems with a single FO term but without time delay. In this paper, a time delayed FO system (in a structure NIOPTD) has been considered and the two different LQR theories for handling time delay systems are applied. This is because the classical LQR theory and the resulting optimal state feedback controller, obtained by solving the Riccati equation are likely to give unstable response in the presence of process delay. The present paper also extends the state-of-the-art techniques by coupling the LQR theory with MOO based determination of weighting matrix ($Q$) and weighting factor ($R$) of the Continuous Algebraic Riccati Equation (CARE) for time delay systems. The optimum weight selection approach is then applied to tune the FOPID gains as the optimal state feedback gains along with integro-differential orders using an incommensurate FO state space formulation, while also keeping the flexibility of choosing FOPID orders independently unlike the approaches reported in [41], [39].

*2.4. Motivation of the present approach*

The motivation of this work is to bridge the gap between the linear quadratic optimal control theory for time delay systems and FO process control using $PI^\lambda D^\mu$ structure to handle a generalized FO template with time delay – NIOPTD [9], [14–18]. This NIOPTD template is capable of capturing the higher order dynamics of a wide variety of self-regulating processes [2] compared to the conventional integer order process models like First Order Plus Time Delay (FOPTD) or Second Order Plus Time Delay (SOPTD) and is also capable of capturing both sluggish and oscillatory dynamics. These compact FO models has the capability of explaining a more generalized power law decaying envelope and Mittag-Leffler oscillation, instead of an exponential envelope and sinusoidal oscillation, commonly encountered in impulse response of integer order ODEs [42]. The present LQR based FOPID controller design first converts the problem in an incommensurate FO state space framework by considering the error signal and the two FO integro-differential orders of the loop error as the state variables. In the present approach, the diagonal elements of the weighting matrix ($Q$) and weighting factor ($R$) are chosen as the decision variables along with the FOPID integro-differential orders ($\lambda$, $\mu$) using the NSGA-II multi-objective optimizer [43]. The optimal state-feedback gains (here the $PI^\lambda D^\mu$ controller gains) associated with the three state variables are then obtained by solving the CARE for each stable solution in the MOO framework. Integral performance indices – ITSE and ISDCO are used here to show that there is a design trade-off between these two conflicting objectives i.e. the set-point tracking performance and required controller effort [15]. These two conflicting objective functions are simultaneously minimized using the multi-objective NSGA-II algorithm to get the Pareto optimal fronts showing the bound of control performance achieved using the two time-delay handling formulation of LQR-FOPID. We also report exhaustive simulation results for oscillatory and sluggish higher order processes with lag-dominant, delay dominant and balanced lag-delay dynamical characteristics in the delayed FO template, using two different delay handling techniques within the LQR framework. This has been translated then to simple FOPID tuning rules for obtaining five controller parameters using the knowledge of the process parameters.

## 3. Theoretical formulation for optimal fractional order controller design

### *3.1. Lyapunov stability to optimal LQR regulator design for fractional order systems*

The classical optimal state-feedback controller or LQR minimizes the infinite horizon quadratic cost function (1) involving the state variables ($x$) and control actions ($u$).

$$J = \int_0^\infty \left[ x^T(t)Qx(t) + u^T(t)Ru(t) \right] dt \qquad (1)$$

Here, $\{Q,R\}$ are the symmetric positive semi-definite weighting matrix and the positive weighting factor respectively, which balance the penalty on the excursion of state variables and control action. Minimization of the integral performance index (1) yields the continuous time algebraic Riccati equation given by (2) which can be used to devise the optimal state-feedback control law (3).

$$A^T P + PA - PBR^{-1}B^T P + Q = 0 \qquad (2)$$

$$u(t) = -R^{-1}B^T Px(t) \qquad (3)$$

where, $\{A,B\}$ are the system matrices of the standard integer order system structure $\dot{x} = Ax + Bu$ and $P$ is the symmetric positive definite solution of the algebraic Riccati equation (2).

It is well-known that minimization of the infinite-horizon cost-function given by (1) leads to the standard Riccati equation in (2) for integer order systems. Now the question arises whether a similar framework holds for fractional order systems or not. Using Lyapunov stability theory to find out the optimal state-feedback control law for integer order systems is an age-old topic in the area of optimal control which involves two steps – (i) formulating the infinite horizon cost function (ii) deriving the ARE using the Lyapunov stability equation. The extension of Lyapunov stability analysis for FO time varying and nonlinear systems has been proposed by Aguila-Camacho *et al.* [44]. Whereas, a different approach of has been adopted by Trigeassou *et al.* [45] which led to a new concepts of the Lyapunov stability analysis for FO systems like sufficiency of quadratic cost function, possible choice of fractional derivative of the energy function [46], generalization for commensurate and incommensurate FO system etc. The Lyapunov stability based Linear Matrix Inequality (LMI) stabilization schemes for FO systems have been discussed in [47], [48].

Sierociuk and Vinagre [38] revisited the earlier works on fractional LQR [37] starting from Euler-Lagrange equation to Riccati differential equation and gave two formulations by considering that the control process is close to the final time ($T_f$) or the half-time ($T_f/2$). By introducing the Riesz potential, the work reported in [38] has shown that the same Riccati equation can be obtained similar to the integer order system with two additional assumptions - (i) the system is around the middle of the control process $t = T_f/2$ and (ii) the time varying matrix $P(t)$ converges to a constant value as $t \to \infty$. Under these assumptions an FO system having structure ${}_0^C D_t^\gamma x = Ax + Bu$ will also produce the same Riccati equation and optimal control law as in (2)-(3). The goal here is to develop a trade-off design using the existing theoretical knowledge of LQR for FO systems [42], [43] which will finally yield optimal FOPID controller parameters to handle FO time delay processes.

Here, the ${}_0^C D_t^\gamma$ denotes the Caputo fractional derivative [7][8] with zero initial condition of order $\gamma$ and is given by (4).

$$ {}_0^C D_t^\gamma f(t) = I^{m-\gamma} D^m f(t) = \frac{1}{\Gamma(m-\gamma)} \int_0^t \frac{D^m f(\tau)}{(t-\tau)^{\gamma-m+1}} d\tau, \quad m = \lceil \gamma \rceil, m \in \mathbb{Z}, \gamma \in \mathbb{R} \quad (4) $$

where, '*C*' stands for Caputo definition, '*I*' stands for integral and '*D*' stands for derivative and *f(t)* represents the function in time domain which undergoes fractional derivative operation.

### *3.2. State-feedback approach of FOPID controller tuning to handle fractional order systems*

Classical PID controller can be designed using the LQR technique, satisfying the quadratic cost function (1) where the state feedback gains can be considered as the PID controller gains [5], [19], [39] if the control loop error, its derivative and integral are considered as the three state variables. Similarly, the concept could be extended as an efficient tuning technique for FOPID controllers using the FO version of LQR [38], if the error and its fractional differ-integrals are considered as the state variables [40]. The formulation of the LQR based FOPID controller for controlling a class of fractional order plant in NIOPTD structure has been shown in Figure 1, where the FO differ-integrated error signals have been considered as the state variables. Here, the class of FO systems is considered to have a NIOPTD template, since most of the higher order oscillatory or sluggish processes can be compactly represented by this structure $Ke^{-Ls}/(Ts^\alpha + 1)$, in terms of the pseudo time-constant, dc gain, time delay and compact system order i.e. $\{T, K, L, \alpha\}$ respectively.

Here, the task is to design an optimal state feedback regulator based FOPID controller that can handle this typical class of FO systems with time delay. In the design of a FOPID controller for delay free systems, the controller gains can be obtained from LQR satisfying the Riccati equation in (2) to produce the control action given by (3). Now the performance of the FOPID control loop could be manipulated in two different ways – (i) the state variables *x(t)* in (3) could be manipulated by the choice of fractional order integral and derivative operators $\{\lambda, \mu\}$, (ii) choice of LQR weights affects the state feedback gains which in turn affects the overall closed loop performance.

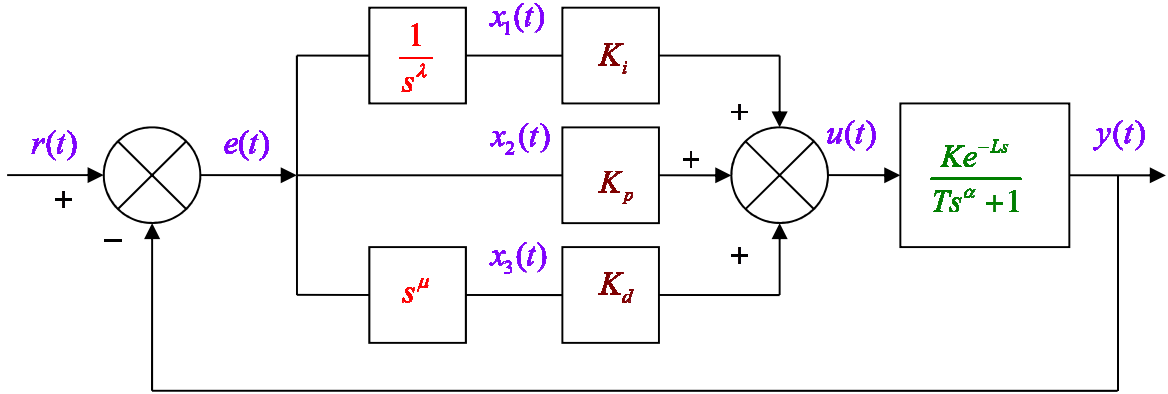

Figure 1: The feedback control system comprising of a NIOPTD plant and PI$^\lambda$D$^\mu$ controller.

In fact, the choice of weighting matrices ($Q$ and $R$) and also the integro-differential orders of fractional state variables of the $PI^\lambda D^\mu$ controller ($\lambda$ and $\mu$) does not affect the optimal regulator formulation given by (1)-(3) for systems with no delay ($L=0$). For each choice of $\{Q,R\}$ it is possible to find out a Riccati solution ($P$) and the associated controller gains $\{K_p, K_i, K_d\}$ by solving the CARE in (2). It is obvious that the closed loop performance changes a lot with variation in all the key parameters $\{Q, R, \lambda, \mu\}$. So, the objective is to find out an optimal set of parameters having the best closed loop tracking performance among all the optimal controllers that could be designed *via* LQR. In order to achieve this, the optimal regulator's weight selection has been improved with other time domain optimality criteria like ITSE and ISDCO (discussed in the next subsection). The ITSE criterion is chosen since it minimizes the overshoot and increases the speed of transient response, while the ISDCO criterion reduces the total variation of the control signal or violent perturbation of the manipulated variable [15].

In Figure 1, if the system is excited with an external input $r(t)$ to have a FO control input $u(t)$ and process output $y(t)$, then the state variables can be considered as (5).

$$x_1(t) = \frac{d^{-\lambda}}{dt^{-\lambda}}[e(t)], \quad x_2(t) = e(t), \quad x_3(t) = \frac{d^\mu}{dt^\mu}[e(t)] \qquad (5)$$

where, differ-integral operator $D = d/dt$ and $e(t)$ is the control loop error. Therefore,

$$\frac{d^\lambda}{dt^\lambda}[x_1(t)] = x_2(t), \quad \frac{d^\mu}{dt^\mu}[x_2(t)] = x_3(t) \qquad (6)$$

It is shown by He *et al.* [5] that in the case of such output feedback control design, the external set-point does not affect the optimal regulator design so that the external set-point change could be considered zero i.e. $r=0$. Thus the closed loop output feedback problem is reduced to a state-feedback regulator problem. Clearly, for set-point $r=0$, the error signal becomes $e = r - y = -y$. Thus, the output signal $y = -e = -x_2$. Now, the FO process in [15][12] with a time delay term is given by (7).

$$\begin{aligned}&\frac{Y(s)}{U(s)} = \frac{Ke^{-Ls}}{Ts^\alpha + 1} \\ \Rightarrow\ &-\frac{E(s)}{U(s)} = \frac{Ke^{-Ls}}{Ts^\alpha + 1} \\ \Rightarrow\ &TD^\alpha[x_2(t)] + x_2(t) = -Ku(t-L)\end{aligned} \qquad (7)$$

Assuming zero initial condition for fractional derivative composition rule, (7) reduces to (8).

$$\begin{aligned}&TD^{\alpha-\mu}x_3(t) + x_2(t) = -Ku(t-L) \\ \Rightarrow\ &D^{\alpha-\mu}x_3(t) = -\frac{1}{T}x_2(t) - \frac{K}{T}u(t-L)\end{aligned} \qquad (8)$$

From (6)-(8) we get (9).

$$\begin{bmatrix} D^{\lambda}x_{1}(t) \\ D^{\mu}x_{2}(t) \\ D^{\alpha-\mu}x_{3}(t) \end{bmatrix} = \begin{bmatrix} 0 & 1 & 0 \\ 0 & 0 & 1 \\ 0 & -1/T & 0 \end{bmatrix} \begin{bmatrix} x_{1}(t) \\ x_{2}(t) \\ x_{3}(t) \end{bmatrix} + \begin{bmatrix} 0 \\ 0 \\ -K/T \end{bmatrix} u(t-L) \quad (9)$$

This is the representation of the above dynamical system with the FOPID controller in a generalized incommensurate FO state-space template, i.e.

$$\frac{d^{q}x(t)}{dt^{q}} = Ax(t) + Bu(t-L) \quad (10)$$

Here, $\frac{d^{q}}{dt^{q}} = \begin{bmatrix} \frac{d^{\lambda}}{dt^{\lambda}} & \frac{d^{\mu}}{dt^{\mu}} & \frac{d^{\alpha-\mu}}{dt^{\alpha-\mu}} \end{bmatrix}^{T}$ and $x(t) = \begin{bmatrix} x_{1}(t) & x_{2}(t) & x_{3}(t) \end{bmatrix}^{T} = \begin{bmatrix} \frac{d^{-\lambda}e(t)}{dt^{-\lambda}} & e(t) & \frac{d^{\mu}e(t)}{dt^{\mu}} \end{bmatrix}^{T}$.

The FO system matrices are given by (11).

$$A = \begin{bmatrix} 0 & 1 & 0 \\ 0 & 0 & 1 \\ 0 & -1/T & 0 \end{bmatrix}; B = \begin{bmatrix} 0 \\ 0 \\ -K/T \end{bmatrix} \quad (11)$$

Since the augmented state space framework (9) has different fractional orders, with a scope of independently choosing the systems and controller orders, the system can be viewed as an incommensurate FO linear dynamical system. If the system is delay free ($L=0$), the fractional state space reduces to the following simple form (12) which is similar to the quadratic optimal control problems often referred in the literatures of analogous integer order systems.

$$\frac{d^{q}x(t)}{dt^{q}} = Ax(t) + Bu(t) \quad (12)$$

Equation (11) represents the open loop system combing both the system and the controller in a generalized incommensurate FO state-space model. The stability criteria for such FO system is not the same as that of the integer order system. Even the right-half plane poles or eigenvalues of the closed loop matrix ($A_{c} = A - BF$) can be stable since the instability region is squeezed from $[-\pi/2, \pi/2]$ to $[-q\pi/2, q\pi/2]$, if it is assumed that state space representation is in commensurate form with order $q$, as a special case of the incommensurate augmented FO system. In addition, the open loop system in (11) should not be judged alone for stability since the LQR formulation with proper weight selection inside the optimization routine calculates the state feedback gains ($F$) which again enforces the stability of the closed loop system using ITSE vs. ISDCO criteria.

It is also important to emphasize that the general solution of the FO systems converges to the classical integer order Riccati equation only under the assumption that the system is around middle of the control process and the time varying Riccati differential equation has a steady-state solution [38]. Along with these two assumptions, let $Q = \begin{bmatrix} Q_{1} & 0 & 0 \\ 0 & Q_{2} & 0 \\ 0 & 0 & Q_{3} \end{bmatrix}$ and $P = \begin{bmatrix} P_{11} & P_{12} & P_{13} \\ P_{12} & P_{22} & P_{23} \\ P_{13} & P_{23} & P_{33} \end{bmatrix}$ which are used to solve the CARE in (2). For a guess value of the weighting matrix $Q$ and $R$, the elements of the positive definite Riccati solution matrix i.e. $\{P_{11}, P_{22}, P_{33}, P_{12}, P_{13}, P_{23}\}$ can be obtained using MATLAB function *lqr()*. Therefore, the state feedback gain matrix ($F$) is obtained as:

$$F = R^{-1}B^T P = \frac{1}{R}\begin{bmatrix} 0 & 0 & -\frac{K}{T} \end{bmatrix}\begin{bmatrix} P_{11} & P_{12} & P_{13} \\ P_{12} & P_{22} & P_{23} \\ P_{13} & P_{23} & P_{33} \end{bmatrix}$$
$$= R^{-1}\begin{bmatrix} \left(-\frac{K}{T}P_{13}\right) & \left(-\frac{K}{T}P_{23}\right) & \left(-\frac{K}{T}P_{33}\right) \end{bmatrix} \quad (13)$$
$$= \begin{bmatrix} -K_i & -K_p & -K_d \end{bmatrix}$$

Since the state variables are chosen in such a way that it represents the error signal and its fractional differ-integrals, the design of the optimal state feedback regulator yields the PI$^\lambda$D$^\mu$ controller gains as the optimal state feedback gain matrix ($F$). The corresponding optimal control law is given by (14).

$$u(t) = -Fx(t) = -R^{-1}B^T Px(t)$$
$$= -\begin{bmatrix} -K_i & -K_p & -K_d \end{bmatrix}\begin{bmatrix} \frac{d^{-\lambda}e(t)}{dt^{-\lambda}} & e(t) & \frac{d^\mu e(t)}{dt^\mu} \end{bmatrix}^T \quad (14)$$
$$= K_p e(t) + K_i \underbrace{\int \cdots \int}_{\lambda-fold} e(t)dt + K_d \frac{d^\mu e(t)}{dt^\mu}$$

### *3.3. Modification of LQR formulation for time delay systems*

It is well known that the LQR framework cannot lead to a guaranteed stabilizing controller if the system has inherent time delay. In standard process control problems, suitable control algorithms like Smith predictors, dead-time compensator, model predictive control etc. are commonly employed for handling large process delays [49]. Optimization based time domain and frequency domain FOPID controller design in industrial process control having large time-delay are discussed in [7][9][14]. Since, the task here is to design an optimal state feedback regulator based FOPID controller, the LQR framework needs to be suitably modified that can stabilize the large process delay without the loss of the optimality condition. We here adopt two such design frameworks which preserves the optimality of state-feedback regulator even for a time-delay system *viz.* (i) by fusing the time delay in the system matrices, proposed by Cai *et al.* [50] and (ii) multiplying nominal delay-free state-feedback gain with an exponential term by He *et al.* [5].

### *3.3.1. Fusion of time delay with system matrices*

Cai *et al.* [50] derived an algorithm to handle systems with time delay using optimum LQR controllers by modifying the system matrices. In this approach, it is shown that the continuous time linear time invariant (LTI) system with delay ($L$) of the structure $\dot{x}(t) = Ax(t) + Bu(t-L)$ can be modified using suitable transformations to produce an LQR framework similar to a delay-free system where the modified system matrices capture the effect of the explicit time delay term. In such a case, the CARE is derived with the augmented system matrices as in (15).

$$A^T \bar{P} + \bar{P}A - \bar{P}[B(A)]R^{-1}[B(A)]^T \bar{P} + Q = 0 \quad (15)$$

where, the modified input matrix for the time delay systems is represented by (16) and the time delay ($L$) appears within the modified input matrix $[B(A)]$ as an exponential term containing the product of stability matrix $A$ and the process delay $L$.

$$B(A) = e^{-AL}B \qquad (16)$$

It is clear that equation (15) is given by the same algebraic Riccati equation (2) except the original input matrix ($B$) is replaced by the modified input matrix $[B(A)]$ to handle the explicit time delay term in the state equation. It has been shown in [50] that in steady state (by neglecting transients of a time varying control term), the control signal is obtained as (17) where $\overline{P}$ is the symmetric positive definite solution of equation (15).

$$u(t) = -\overline{F}x(t) = -R^{-1}[B(A)]^T \overline{P}x(t) = -R^{-1}B^T[e^{-AL}]^T \overline{P}x(t) \qquad (17)$$

In addition, the weighting matrices for the LQR i.e. $\{Q, R\}$ and integro-differential orders $\{\lambda, \mu\}$ are chosen using multi-objective NSGA-II algorithm while minimizing the ITSE representing good set-point tracking and ISDCO representing minimum variation of the manipulated variable to obtain the trade-off between them. With the assumptions of the process being in the middle of the final control time and with a steady state Riccati solution $\overline{P}$, the delayed FO model in (10) can now be reduced to the standard delay free FO state-space (12) with modified input matrix (16). The optimal state-feedback controllers can now be obtained using the standard Riccati equation but with the modified input matrix $[B(A)]$, as shown in (15).

### 3.3.2. Multiplying nominal state-feedback gains with an exponential term containing the product of closed loop matrix and delay

He *et al.* [5] developed a method where the controller gains are initially time varying. Finally as the transient response crosses the system delay ($L$) after a set-point change, the controller gains become constant. The steady state controller gains that stabilize a similar state space model with time delay, has been derived from the LQR theory which results in an additional exponential term along with the traditional state feedback gains $F$ in (13). The steady state control signal for He's method [5] of handling time delay systems with LQR formulation and steady state value of the time varying PID controller gains are given by (18).

$$\begin{aligned} u(t) &= -R^{-1}B^T P e^{(A-BR^{-1}B^T P)L} x(t) \quad , t \geq L \\ &= -F e^{(A-BF)L} x(t) = -F e^{A_c L} x(t) \end{aligned} \qquad (18)$$

Here, the closed loop matrix $A_c = A - BF$, where $F$ refers to the optimal state feedback gains for the delay free system in (12). Therefore, a second formulation is obtained to handle the time delay term within the LQR design. The similarity between these two formulations is that both of them ignores the transient dynamics of the controller which is implemented as a small time varying term within the interval $t < L$ and both of them contains an exponential term involving the time delay of the process. The difference between the control actions in (17) and (18) is that the Riccati solution of Cai's method [50] is based on the modified CARE (15) i.e. $\overline{P}$ involving modified input matrix $[B(A)]$ but in He's method [5] the same Riccati solution (*P*) and state feedback gains (*F*) are used along with some additional terms. This paper compares which formulation among these two is capable of producing a non-dominated Pareto front indicating better control performance.

### *3.4. Optimum selection of weighting matrices of LQR using two conflicting time-domain control objectives*

The above two LQR formulations are optimal for a specific choice of the weighting matrices *Q* and *R*. Indeed, the time domain performance is heavily affected for any arbitrary choice of the LQR weights, although in each case the state optimality is preserved from (1). This is logical since the choice of weighting matrices determine the state feedback gains ($PI^\lambda D^\mu$ controller gains in this case) that directly affect the performance of the closed loop system. In order to handle this problem, a MOO technique is employed by minimizing two conflicting time domain performance indices $J_1$ and $J_2$ (19) as also studied in [15]. This tunes the elements of the diagonal elements of the weighting matrix and weighting factor i.e. $\{Q_1, Q_2, Q_3, R\}$ [19] and integro-differential orders of the FOPID controller i.e. $\{\lambda, \mu\}$ [40].

$$J_1^{ITSE} = \int_0^\infty te^2(t)dt, \quad J_2^{ISDCO} = \int_0^\infty \left(u(t) - u_{ss}\right)^2 dt \quad (19)$$

Such a choice of control objectives has been reported to have a design trade-off since higher speed of tracking for the process variable needs large perturbation of the manipulated variable or the control signal [15], [14]. Thus instead of focusing on a particular controller minimizing a single objective function as weighted sum of two objectives ITSE and ISDCO [14], a set of controllers that lie on the Pareto front is sought, since optimal choice of such weights of the two parts of the cost function are often heuristic and largely depends on the controller and the system under control [15]. Thus some of the solutions may give good performance in terms of one objective at the cost of deterioration in the others, while for other solutions on the Pareto front it would be vice-versa. Hence the set of non-dominated solutions obtained on the Pareto front give the limits of the controller performance. Thus the controller solutions lying on the Pareto front cannot perform better in one control objective without a corresponding deterioration in the performance of the other.

The rationale for using these specific integral performance indices in (19) is to get a good time domain response and at the same time to limit the deviation in the controller output to avoid actuator saturation and integral wind-up [3]. Instead of the ITSE criterion for set-point tracking, other criterion like Integral of Time multiplied Absolute Error (ITAE) or Integral of Absolute Error (IAE) could have been used which would have resulted in smaller penalty for high oscillations at later stages [9][2]. The deviation of the control signal is also minimized in the form a performance index known as ISDCO to limit violent perturbation of the manipulated variable. At a first glance this might seem as a redundant criteria since the LQR methodology already gives optimal values of the controller gains with the lowest cost. However, this is actually obtained for a fixed value of the weighting matrices. When *Q* and *R* are varied, for each choice of weighting matrices, the LQR would give an optimal gain with the lowest possible cost, but that does not necessarily imply a good time domain performance [4], [19], [21], [22] with the LQR cost function (1). Also, for an optimal choice of weighting matrices (*Q* and *R*) and differ-integral orders ($\lambda$ and $\mu$), the FOPID tuning problem becomes optimal due to the introduction of time domain performance indices (19) as well as the classical optimal regulator (LQR) based approach (1), involving the fractional states.

*3.5. Multi-objective optimization framework for tuning LQR weights and FOPID orders*

A generalized multi-objective optimization framework can be defined as follows:

$$\text{Minimize } F(x) = (f_1(x), f_2(x), ..., f_m(x))$$
$$\text{subject to: } g_i(x) \leq 0, \forall i \in [1, p], \quad (20)$$
$$h_j(x) = 0, \forall j \in [1, q]$$

such that $x \in \Omega$.

where, $\Omega$ is the decision space, $\mathbb{R}^m$ is the objective space, $F: \Omega \rightarrow \mathbb{R}^m$ consists of $m$ real valued objective functions and $g_i(\cdot)$ and $h_j(\cdot)$ are the optional $p$ number of inequality and $q$ number of equality constraints on the problem respectively.

Let, $u = \{u_1, ..., u_m\} \in \mathbb{R}^m$ and $v = \{v_1, ..., v_m\} \in \mathbb{R}^m$ be two vectors and $u$ is said to dominate $v$ if $u_i < v_i \forall i \in \{1, 2, ..., m\}$ and $u \neq v$. A point $x^* \in \Omega$ is called Pareto optimal if $\nexists x | x \in \Omega$ such that $F(x)$ dominates $F(x^*)$. The set of all Pareto optimal points, denoted by PS is called the Pareto set. The set of all Pareto objective vectors, $PF = \{F(x) \in \mathbb{R}^m, x \in PS\}$, is called the Pareto Front. This implies that no other feasible objective vector exists which can improve one objective function without simultaneous worsening some other objective function.

The NSGA-II algorithm [43] converts $m$ diverse objectives into one single fitness function by creating a number of different fronts. The solutions on these fronts are refined iteratively based on their distance with their neighbours (crowding distance) and their level of non-domination. The NSGA-II algorithm ensures that the solutions found are close to the original Pareto front and are diverse enough to find the whole length of the Pareto front.

For the present simulation study, the population size is taken as 100 and the number of generations as 100. The elite count represents the number of fittest individuals which are directly copied over to the next generation. The Pareto fraction in this study is considered as 0.7. An intermediate crossover scheme is adopted which produces off-springs by random weighted average of the parents. The mutation scheme adds a random number at an arbitrary point in the individual. The variables that constitute the search space for the LQR based FOPID design are $\{Q_1, Q_2, Q_3, R, \lambda, \mu\}$. The intervals of the search space for these variables are $\{Q_1, Q_2, Q_3, R\} \in [0, 100]$ and $\{\lambda, \mu\} \in [0, 2]$. In fact, with this search interval, the number of state variables remains always three in the state space formulation (9) even though the integral and differential orders may take values higher than unity, representing faster time response and better loop stability respectively. The MOO algorithm is terminated when the average change in the Pareto spread of the generation is not significant. It is also possible to encounter local minima in the objective function space within the MOO. In order to ensure that a true global minima has been found, we ran the algorithm multiple times and report here the best result with the most non-dominated Pareto front. Also, in the present design scenario, the LQR weights does not need any particular initialization as they are chosen with an MOO algorithm which automatically takes random initial guess values every time. The initialization of LQR weights is an important issue when the optimization is a deterministic one like any gradient-descent algorithm and the application is intended for online implementation for example in [51].

Also unstable response within the optimization routine shows high values of the chosen cost functions – ITSE and ISDCO and are therefore automatically rejected if a bad random guess is encountered. This way the MOO evolves over the iterations to yield better solutions by simultaneously minimizing both the design objectives – finally giving the optimal Pareto front on which all the solutions are not only stable but also represent the best possible trade-off between the two conflicting objectives.

## 4. Illustrative examples

The FO plants that have been considered here, show heavily oscillatory and sluggish open loop response as shown in Figure 2. The process parameters of the two test plants in NIOPTD structure have been chosen from the study by Ruszewski [52]. Figure 2 also shows that the FO processes (7) exhibit sluggish and oscillatory open loop dynamics for $(\alpha < 1)$ and $(\alpha > 1)$ respectively which is also evident from the sharp increase in the system $H_\infty$ norm for $(\alpha > 1)$.

It is observed that the system shows oscillatory ($2 > \alpha > 1$) and sluggish ($1 > \alpha > 0$) open loop response, even with the simple first order plus time delay (FOPTD) like FO template due to the presence of higher order dynamics of the plant which can be easily and compactly modelled as a FO transfer function with delay. The present simulation studies are reported using the FOMCON toolbox [53] and the performance measures are calculated on a finite time horizon of 100 sec.

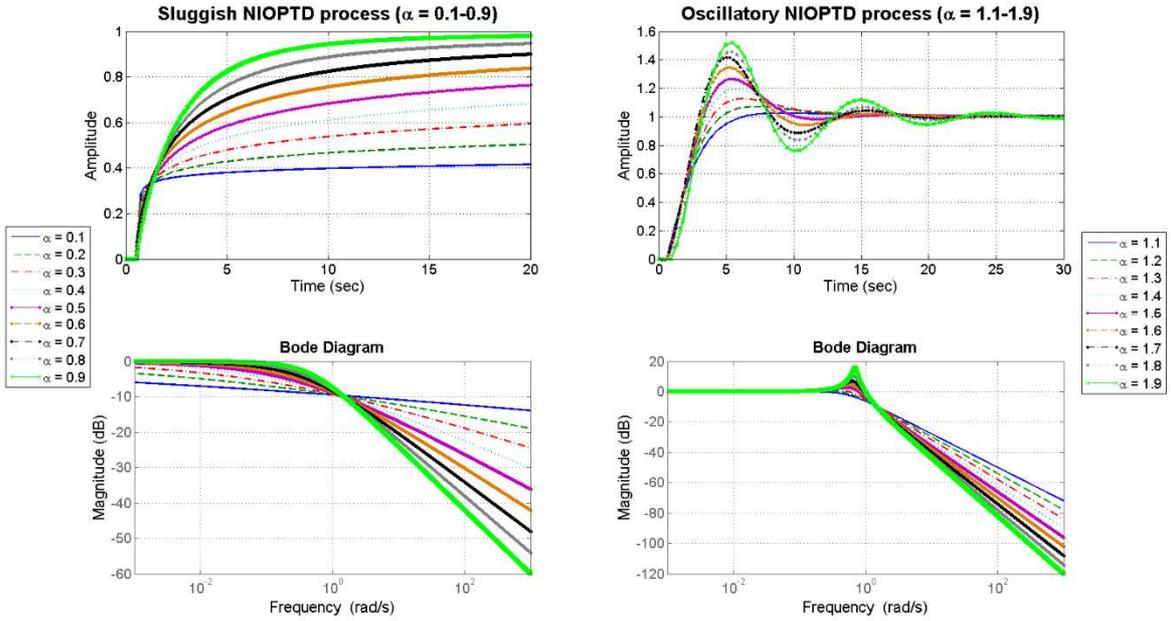

Figure 2: Open loop step response and Bode magnitude diagram of NIOPTD plant (7) with *K=1*, *L=0.5*, *T=2* with change in *α*, exhibiting sluggish and oscillatory dynamics.

### 4.1. Oscillatory fractional order process with time delay

The oscillatory system under consideration is represented by the following lag-dominant ($T > L$) transfer function [52]:

$$G_1(s) = \frac{1e^{-0.5s}}{2s^{1.5} + 1} \tag{21}$$

Figure 3 shows the obtained Pareto fronts for the NIOPTD process using both Cai's method [50] and He's method [5]. It is observed that both the methods can handle time delays and give a set of non-dominated solutions on the Pareto front. For the set of solutions obtained with Cai's method, lower values of ITSE are obtained at the cost of higher values of ISDCO. On the other hand, for He's method, the Pareto optimal solution set contains lower values of ISDCO at the expense of higher values of ITSE. However, in He's method, the solutions lie on the convex side of those with Cai's method. This implies that the solution obtained using He's method are better (non-dominated with respect to Cai's method). Next, three representative solutions are chosen on the Pareto front with He's method i.e. $\{A_1, B_1, C_1\}$ and the corresponding time domain simulations are shown as well. These solutions are the ones at the extreme end and the median solution on the Pareto front.

Figure 4 show the time domain evolution of the process output, controller output and states for the three representative solutions on the Pareto front $\{A_1, B_1, C_1\}$ using the non-dominated He's method of LQR-FOPID. A load disturbance is applied at the later stages of settling down the oscillation for set-point change and the disturbance rejection properties of the obtained controller are also investigated. As can be inferred from the Pareto front and also the time domain simulations, the solution $C_1$ has the highest overshoot and settling time than the solution $A_1$, whereas the solution $B_1$ lies in between. However it is the other way round for the control signal. Solution $A_1$ has the highest deviation in control signal and $C_1$ has the lowest. Hence the simulation results verifies our proposition that the set point tracking and the control signal are conflicting objectives. Solution $A_1$ gives the best load disturbance rejection over $B_1$ and $C_1$. However, the load disturbance rejection was not explicitly taken into the optimization framework unlike [14] since its integration within an LQR framework is not very popular. The simulation studies show that the solutions would work in a practical setting as well, as physical processes must be able to reject load disturbances to a sufficient level for effective functioning. Since the LQR based method also minimizes the deviation in the state trajectories, the state variables i.e. the loop error and their fractional differ-integrals have also been shown in Figure 4.

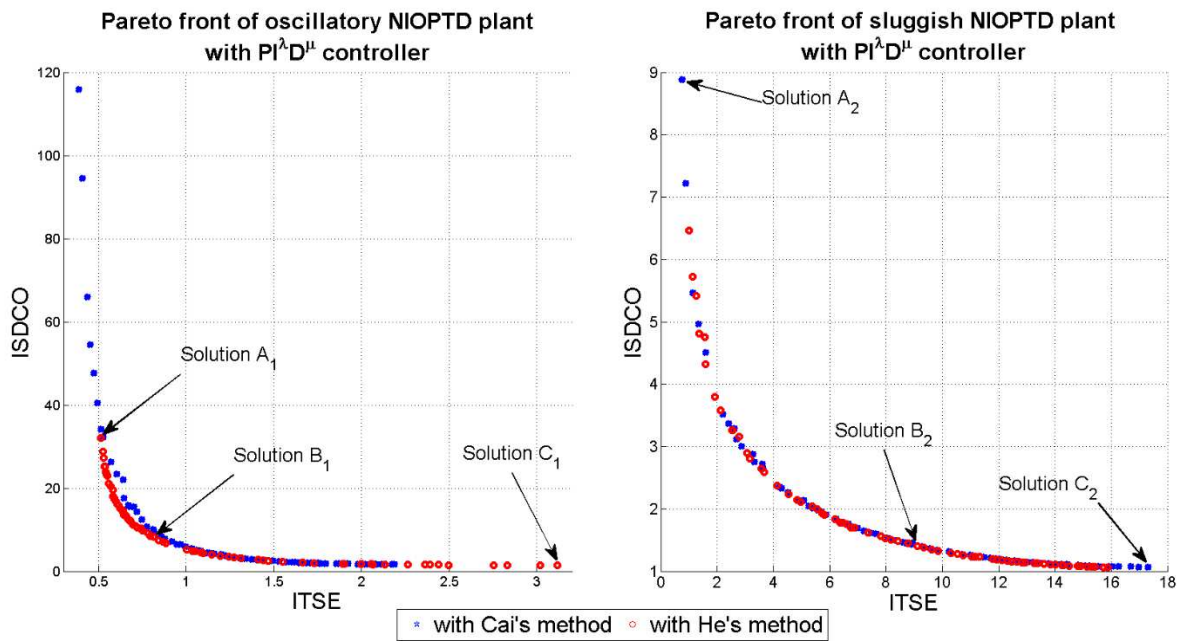

Figure 3: Pareto fronts of NIOPTD plant with FOPID controller using two methods of handling time delay in LQR: (left) oscillatory process (21) (right) sluggish process (22).

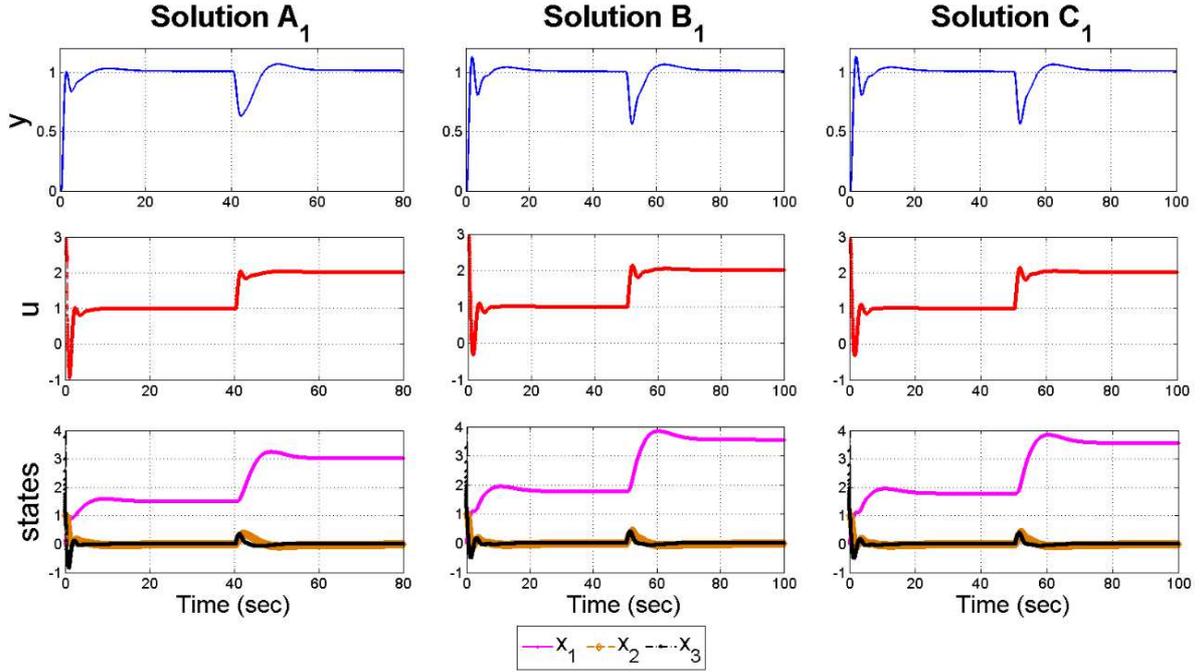

Figure 4: Time response, controller output and state trajectories of oscillatory NIOPTD plant (21) with LQR based FOPID controller for solution $A_1$

*4.2. Sluggish fractional order process with time delay*

A sluggish process is considered next having the lag-dominant ($T > L$) transfer function (22), similar to that studied in [52].

$$G_2(s) = \frac{1e^{-0.5s}}{2s^{0.5} + 1} \qquad (22)$$

Figure 3 shows the Pareto fronts for the sluggish NIOPTD processes with both the Cai's method [50] and He's method [5] for handling the time delay in an LQR framework. As can be seen in this case, the Pareto front obtained using Cai's method has a larger spread, i.e. has got more diverse non-dominated LQR-FOPID solutions. Therefore, three representative solutions ($A_2$, $B_2$, $C_2$) are chosen on this Pareto front and the time domain responses are plotted for each of the cases in Figure 3 respectively. As can be seen from the Pareto front and also from the time domain simulations, the Solution $A_2$ has the fastest settling time (within 20 seconds), while it is longer for solutions $B_2$ (settling time within 40 seconds) and solution $C_2$ (settling time within 100 seconds). However an investigation into the time domain evolution of the control signals would show that solution $A_2$ has a much higher control signal than solutions $B_2$ and $C_2$. This again reaffirms the proposition that better time domain performance can only be obtained at the expense of increased control cost [14]. The load disturbance rejection performances for the three solutions show that $A_2$ performs the best in this category, followed by $B_2$ and lastly $C_2$.

Table 1 and Table 2 report the numerical values of the representative solutions on the Pareto fronts which are obtained after multi-objective optimization. Table 1 shows the three representative solutions as the LQR weighting matrices and FOPID integro-differential orders from each of the best Pareto fronts for time delay handling i.e. using He's and Cai's method to control the oscillatory and sluggish NIOPTD plants respectively. From the LQR weights given in Table 1, the FOPID gains can be calculated by solving the matrix Riccati equations, corresponding to two different methods discussed for time delay handling and have been reported in Table 2.

Table 1: Representative solutions on the Pareto front showing the LQR matrices for lag-dominant processes using different delay handling methods

| Type of Process | Time delay handling in LQR formulation | Solution | $Q_1$ | $Q_2$ | $Q_3$ | $R$ | $\lambda$ | $\mu$ |
|---|---|---|---|---|---|---|---|---|
| Oscillatory lag dominant ($\alpha = 1.5$) | He's method [5] | $A_1$ | 0.970396 | 0.040181 | 0.022387 | 0.204583 | 1.071069 | 0.716467 |
| | | $B_1$ | 0.643793 | 0.02965 | 0.062444 | 0.34342 | 1.133782 | 0.449655 |
| | | $C_1$ | 0.086837 | 0.023281 | 0.095594 | 0.992322 | 1.382362 | 0.035294 |
| Sluggish lag dominant ($\alpha = 0.5$) | Cai's method [50] | $A_2$ | 0.605858 | 0.080236 | 0.057087 | 0.946696 | 0.995725 | 0.026867 |
| | | $B_2$ | 0.061832 | 0.033902 | 0.09303 | 0.873642 | 0.891239 | 0.026349 |
| | | $C_2$ | 0.049785 | 0.026213 | 0.098279 | 0.918109 | 0.754981 | 0.026134 |

From the above simulations, apparently it may appear that He's method of LQR-FOPID for oscillatory NIOPTD plants and Cai's method of LQR-FOPID for sluggish NIOPTD plants yields better non-dominated Pareto fronts using ITSE and ISDCO criteria. However, consistent winning of one method over the other for a class of FO process needs to be verified for varying level of dominance between the time delay and time constant (lag), along with variation in the characteristics exponent ($\alpha$) of the FO process which are reported in the following sub-sections.

Table 2: Representative solutions of tuned controller parameters on the Pareto front along with objective function values for the lag-dominant processes

| Type of Process | Solution on the non-dominated Pareto front | ITSE | ISDCO | $K_p$ | $K_i$ | $K_d$ | $\lambda$ | $\mu$ |
|---|---|---|---|---|---|---|---|---|
| Oscillatory lag dominant | $A_1$ | 0.515799 | 32.10448 | 0.6718 | 0.9327 | 2.053 | 1.071069 | 0.716467 |
| | $B_1$ | 0.816633 | 8.217709 | 0.5692 | 0.7092 | 1.8411 | 1.133782 | 0.449655 |
| | $C_1$ | 3.116587 | 1.434095 | 0.2182 | 0.0903 | 0.9434 | 1.382362 | 0.035294 |
| Sluggish lag dominant | $A_2$ | 0.772218 | 8.874867 | 0.8 | 0.9186 | 2.6498 | 0.995725 | 0.026867 |
| | $B_2$ | 8.720682 | 1.452479 | 0.266 | 0.119 | 1.1945 | 0.891239 | 0.026349 |
| | $C_2$ | 17.32365 | 1.067778 | 0.2329 | 0.0791 | 1.0728 | 0.754981 | 0.026134 |

Now, the robustness property of the median solutions of the FOPID controller (with the parameters give in Table 2) on the Pareto front are explored in Figure 6, for the two chosen plants with oscillatory ($B_1$) and sluggish ($B_2$) open loop dynamics. The purpose here is to show how the closed loop performance changes with variation in the time constant (*T*) and delay (*L*), while the FOPID controller parameters are kept fixed. The parametric robustness in Figure 6 is shown in terms of increase in both the ITSE and ISDCO measures for the best delay handling scheme reported in Table 1 and Table 2. As expected, for the oscillatory process ($\alpha = 1.5$), both the performance measures deteriorates much faster than that of the sluggish process ($\alpha = 0.5$). The decrease in the closed loop performance is inevitable while considering a significant deviation from the nominal process parameters which were used for tuning the controller. Figure 6 shows that the ITSE and ISDCO measures does not blow up and is still capable of stabilizing the process as well as keeping the deterioration of closed loop performances within allowable limits.

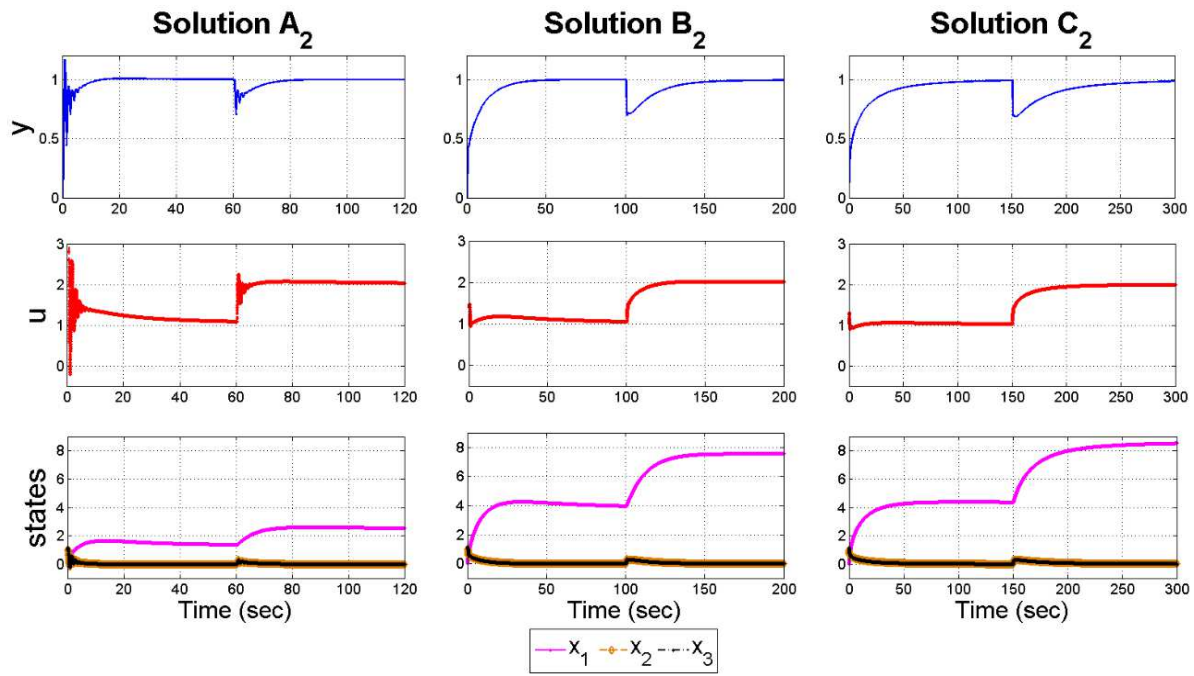

Figure 5: Time response, controller output and state trajectories of sluggish NIOPTD plant with LQR based FOPID controller for solution $A_2$

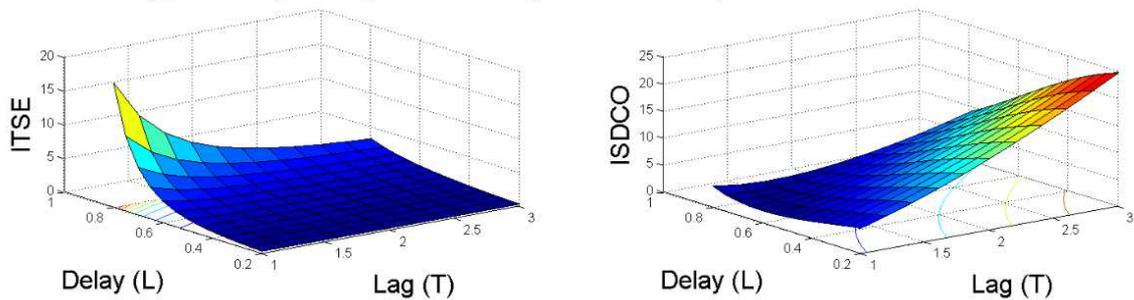

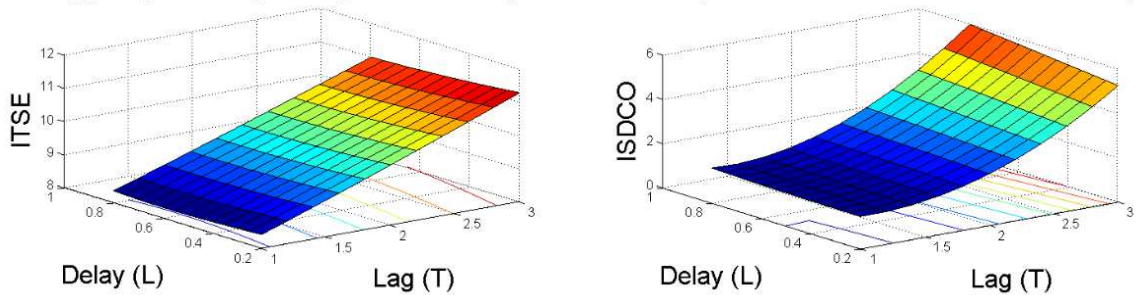

Figure 6: Robustness of the Pareto median solution against variation in process delay and lag.

*4.3. Multi-objective control of lag dominant, balanced lag delay and delay dominant plants*

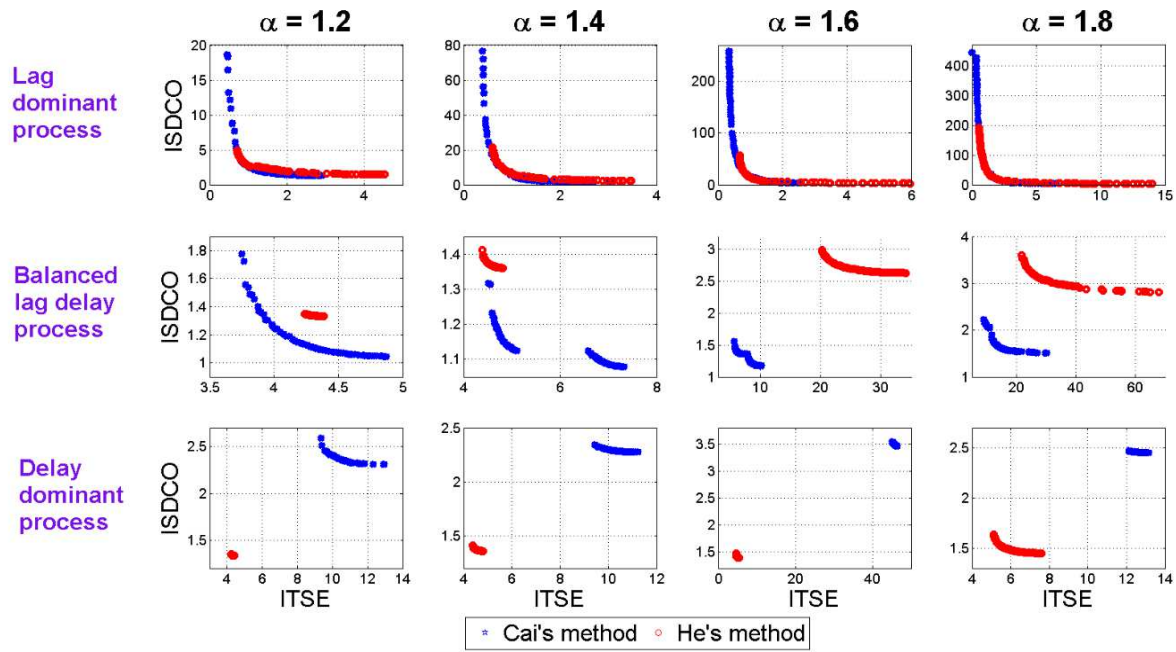

Figure 7: Comparison of the Pareto fronts for the oscillatory NIOPTD processes ($\alpha>1$).

Figure 7 shows the comparison of the Pareto fronts with the two delay handling techniques (Cai's and He's method) for different values of $\alpha$ that determines the oscillatory and sluggish behavior of the open loop system as shown in Figure 2. We here explore which method has yielded a non-dominated Pareto front for the three different class of NIOPTD processes *viz.* lag-dominant ($T>L$), balanced lag and delay ($T \approx L$) and delay dominant ($T<L$) processes. From Figure 7, it is evident that for the balanced lag-delay processes the Cai's method and for the delay dominant processes the He's method gives a better non-dominant Pareto front, although the length of the Pareto front may vary, depending on the oscillatory nature ($\alpha$) of the open loop process. Whereas for the lag-dominant processes the Pareto fronts are comparable and as such there is no clear winner. It is evident that for lag-dominant processes Cai's method keeps the ITSE low but cannot restrict an increase in ISDCO whereas the observation is just the reverse with He's method resulting in low ISDCO but increased ITSE. In some cases, two Pareto front cuts each other which indicates a weak Pareto dominance, as such one controller is better only in a particular regime.

A mixed response is observed for the sluggish NIOPTD processes in Figure 8, especially in the balanced lag-delay and delay dominant processes where the He's method outperforms the Cai's method in providing non-dominated Pareto fronts. But for $\alpha = 0.6$ and $\alpha = 0.8$, Cai's method is slightly better. For these two values of $\alpha$, the lag-dominant processes show comparable results with both the methods. For highly sluggish lag-dominant NIOPTD processes ($\alpha = 0.2$), He's method keeps the ISDCO low and Cai's method yields low ITSE, whereas for moderately sluggish ($\alpha = 0.4$) lag-dominant process, He's method clearly outperforms the other. Also, in both Figure 7 and Figure 8, some of the Pareto fronts are discontinuous indicating forbidden regimes in the MOO space. This discontinuity in the Pareto fronts could be a result of partitioning in the dynamical characteristics for FO elements below and above one, using the rational approximation, or could be an effect of the process itself which does not yield any stabilizing FOPID controller to operate in that particular regime. These discontinuity in the Pareto fronts needs to be investigated in future studies.

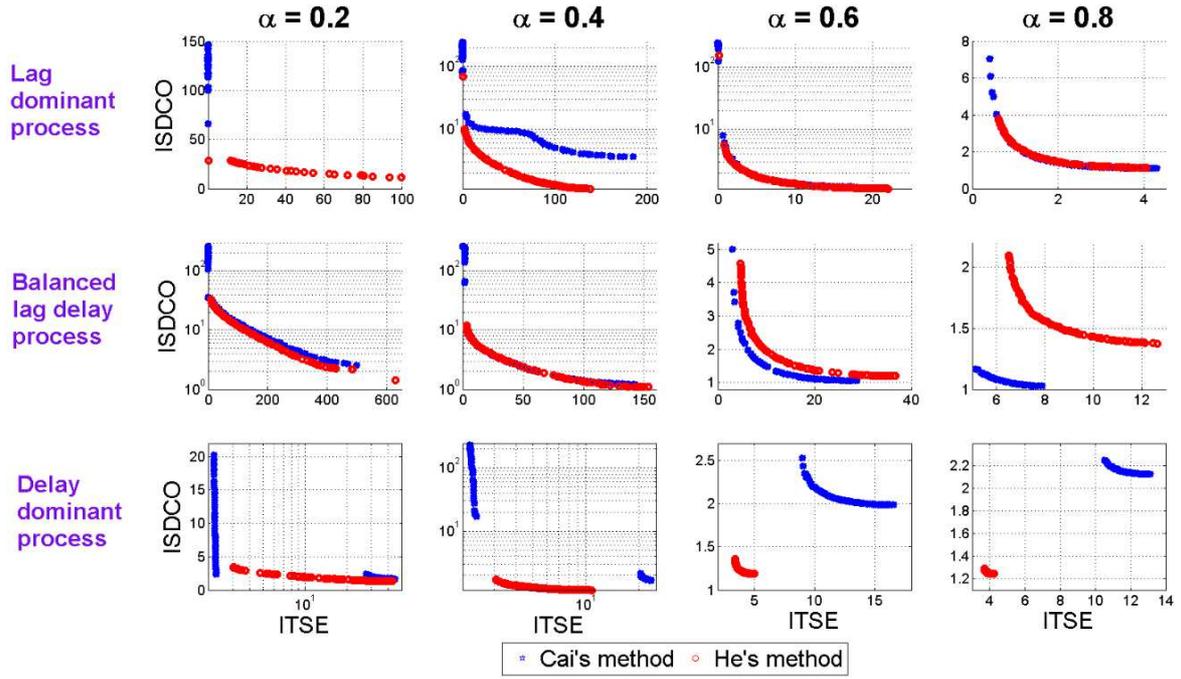

Figure 8: Comparison of the Pareto fronts for the sluggish NIOPTD processes ($α<1$).

Table 3: Summary of the non-dominated Pareto fronts for different NIOPTD process order and dominant characteristics

| NIOPTD process order | Process characteristics | | |
|---|---|---|---|
| | Lag dominant | Balanced lag and delay | Delay dominant |
| $α = 0.2$ | He's method | He's method | He's method |
| $α = 0.4$ | He's method | weak dominance | He's method |
| $α = 0.6$ | He's method | Cai's method | He's method |
| $α = 0.8$ | weak dominance | Cai's method | He's method |
| $α = 1.2$ | weak dominance | Cai's method | He's method |
| $α = 1.4$ | weak dominance | Cai's method | He's method |
| $α = 1.6$ | weak dominance | Cai's method | He's method |
| $α = 1.8$ | weak dominance | Cai's method | He's method |

Comparing the results obtained in section 4.1-4.3, we can conclude that the non-dominated Pareto front between the two delay handling methods not only depend on the process characteristics but also on the nature of its open loop dynamics, characterized by its order $α$. As shown in Table 3, for delay dominant plants, He's method is consistently better. For balanced lag-delay processes, in most of the cases Cai's method is better. For lag-dominant processes, both the methods give good result and are comparable to each other which produces a weak dominance between the two methods.

Especially in industrial control with very high order process dynamics (modelled by the NIOPTD process) and large time delay, the results reported in Figure 8 could be useful in selecting the delay handling methodology with an optimal FOPID controller. It has been shown in [9], [14] that very high order process models can be compactly represented as FO transfer function template NIOPTD. The model reduction results has been established on the higher-order test bench plants, available in contemporary literature [14][9]. Therefore, our results are based on the assumption that the higher order process dynamics could be faithfully represented in the NIOPTD template and we here focus on the results of LQR based FOPID controller to handle such FO processes.

### *4.4. Summary of the proposed approach and translating LQR trade-off design in FOPID tuning rules for NIOPTD plants*

The steps of the proposed LQR-FOPID design algorithm is summarized below:

***Step 1***: Reduce higher order oscillatory/sluggish process dynamics in NIOPTD template (7).

***Step 2***: Run MOO algorithm to select LQR weights $\{Q, R\}$ and FOPID orders $\{\lambda, \mu\}$ using ITSE and ISDCO as the two conflicting objectives (19).

***Step 3***: Select the transformation by Cai's method (17) and He's method (18) for time delay handling within LQR.

***Step 4***: Solve the associated CARE within the MOO to calculate FOPID gains $\{K_p, K_i, K_d\}$ and obtain the Pareto front.

***Step 5***: Select the non-dominated Pareto front between Cai's and He's method.

***Step 6***: Select the median solution on the non-dominated Pareto front as a trade-off design.

***Step 7***: Check for the robustness of the FOPID controller settings by varying the process parameters.

Also, similar to the rule based PID gains selection, recent studies have also suggested using tuning rules to select integro-differential orders $\{\lambda, \mu\}$ of the FOPID controller. Most results have been reported for integer order plants [11], [10], but the counterpart for FO processes have not been investigated yet. However none of the above tuning rules consider the multi-objective LQR formalism and report design trade-offs between conflicting objectives. Also, any such tuning rule to select either the FOPID gains or orders will not give the design trade-offs between the two chosen objectives which is the prime focus of this paper. Rather a specific tuning rule will give a single point on the chosen objective function space. Therefore, the integro-differential orders are selected along with the LQR weights using an MOO algorithm which evolves over the generations yielding a non-dominated Pareto optimal front that can be considered as the achievable optimal trade-off design [15], [54], using a particular controller structure.

Next, the analytical expressions of the five FOPID parameters are to be formulated (including the integro-differential operators) as functions of the FO process parameters $\{L, T, \alpha\}$ – showing both sluggish and oscillatory open loop dynamics. The tuning rules are useful for easy calculation of the FOPID parameters to control the FO plant without running the optimization to get the trade-offs. The tuning rules are also optimum since it balances both the conflicting objectives, as their median solution on the Pareto front. On the other hand, the LQR optimality of the FOPID controller is already enforced within the problem formulation while also efficiently handling the time delays. The median values are selected from the non-dominated Pareto front (between He's and Cai's method) and are reported in the supplementary material. Using the data reported in the supplementary material, the tuning rules are now generated for the FOPID controller tuning knobs, as a function of two significant NIOPTD process parameters – delay to lag ratio ($L/T$) and order of the process ($\alpha$). The goodness of fit measures are also computed for the tuning rules as a test of true representation of the data used to construct the rules as shown in

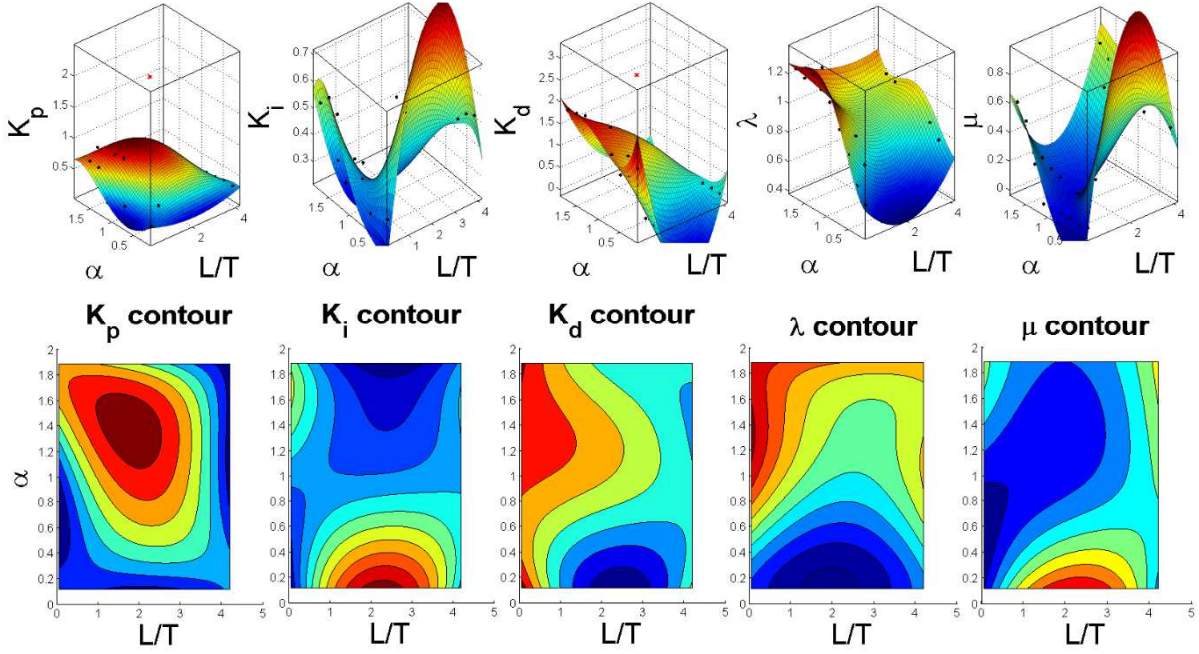

Figure 9: Fitted tuning rule surface and contours for the five FOPID parameters.

. The fitted surface and contour plots of the tuning rules for five FOPID parameters are depicted in Figure 9. One outlier in $K_p$ and $K_d$ are removed before fitting the tuning rule using polynomial regression. In the tuning rule generation process, the orders of both ($L/T$) and $\alpha$ have been varied from one to four and the best model was selected with highest adjusted coefficient of determination (*adjusted $R^2$*). Here the *adjusted $R^2$* is chosen as the deciding statistical measure for model selection since it penalizes more complex models unlike $R^2$, thus reducing the chance of over-fitting. The *adjusted $R^2$* also does not capture the effect of different scales of the independent variables similar to the Root Mean Squared Error (RMSE). Therefore, although more complex models could have been yielded a better fit and push the $R^2$ near one, the *adjusted $R^2$* provides a safe-guard against complex rule generation (using only $R^2$) [11]. It also avoids the discrepancy due to difference in scaling of the two input variables ($L/T$) and $\alpha$ (using only RMSE). The tuning rules for the median non-dominated LQR-FOPID controllers are given in (23) as a polynomial model in ($L/T$) and $\alpha$, of the order of 2 and 4 respectively and the associated coefficients of the polynomial models for all these five parameters are reported in Table 5.

$$\begin{aligned}
\{K_p, K_i, K_d\} &= \left[f(L/T, \alpha)\right]/K, \quad \{\lambda, \mu\} = f(L/T, \alpha) \\
f(L/T, \alpha) &= p_{00} + p_{10}(L/T) + p_{01}\alpha + p_{20}(L/T)^2 + p_{11}(L/T)\alpha + p_{02}\alpha^2 \\
&+ p_{21}(L/T)^2\alpha + p_{12}(L/T)\alpha^2 + p_{03}\alpha^3 + p_{22}(L/T)^2\alpha^2 + p_{13}(L/T)\alpha^3 + p_{04}\alpha^4
\end{aligned} \quad (23)$$

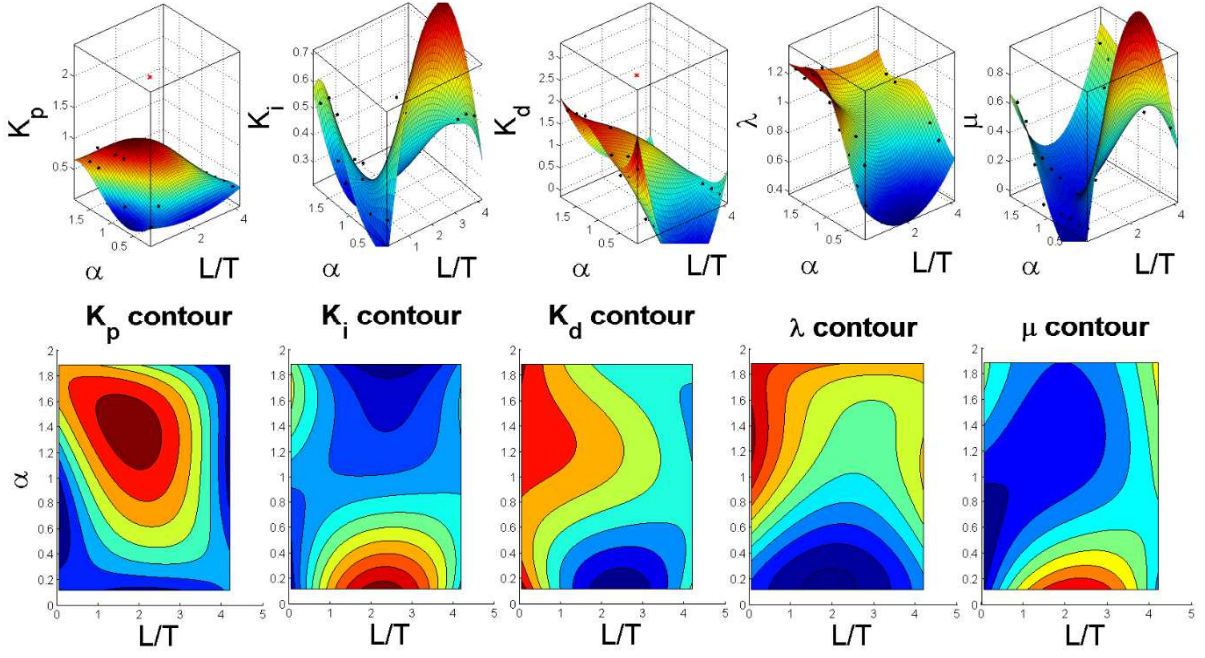

Figure 9: Fitted tuning rule surface and contours for the five FOPID parameters.

## *4.5. Discussions on achievements, assumptions and limitations*

The contribution of this paper, can be considered as the successful integration of LQR theory for the design of optimal FOPID controllers for FO processes with time delay within an MOO framework. As in many other design problems of LQR, the system has been considered to be linear and noise free. The main emphasis is to show that an LQR based FOPID controller could be designed for a class of FO systems – NIOPTD, exhibiting both oscillatory and sluggish open loop dynamics. The design methodology would help to enjoy the state optimality of LQR technique even in the presence of process delay along with the externally imposed Pareto optimality between ITSE and ISDCO by obtaining the best set of trade-off values for the LQR weights and FOPID differ-integral orders.

Table 4: Goodness of fit for the regression model in the tuning rules of five FOPID parameters

| FOPID parameter | Order of ($L/T$) | Order of $\alpha$ | Adjusted $R^2$ | $R^2$ | RMSE |
|---|---|---|---|---|---|
| $K_p$ | 2 | 4 | 0.7875 | 0.8937 | 0.104 |
| $K_i$ | 2 | 4 | 0.7328 | 0.8606 | 0.06512 |
| $K_d$ | 2 | 4 | 0.9837 | 0.9919 | 0.09768 |
| $\lambda$ | 2 | 4 | 0.9535 | 0.9758 | 0.05212 |
| $\mu$ | 2 | 4 | 0.899 | 0.8065 | 0.1268 |

Table 5: Coefficients of the regression model (23) for the tuning rules of five FOPID parameters

| Coefficients | $K_p$ | $K_i$ | $K_d$ | $\lambda$ | $\mu$ |
|---|---|---|---|---|---|
| $p_{00}$ | 0.4225 | 0.001375 | 3.39 | 0.5972 | 0.06535 |
| $p_{10}$ | -0.3738 | 1.002 | -3.976 | -0.1805 | 0.1732 |
| $p_{01}$ | -0.8846 | 0.7251 | -8.749 | -0.3615 | -0.2331 |
| $p_{20}$ | 0.08037 | -0.2251 | 0.8184 | 0.04781 | 0.1506 |
| $p_{11}$ | 2.079 | -1.216 | 7.177 | -0.3342 | -0.2898 |
| $p_{02}$ | 0.05753 | -1.36 | 12.95 | 2.808 | 0.3122 |
| $p_{21}$ | -0.4099 | 0.3161 | -1.484 | 0.08372 | 0.3712 |
| $p_{12}$ | -1.245 | 0.09725 | -3.758 | 0.03983 | 0.01343 |
| $p_{03}$ | 0.935 | 1.389 | -7.427 | -2.304 | -0.05011 |
| $p_{22}$ | 0.1884 | -0.07726 | 0.6184 | -0.05261 | -0.1479 |
| $p_{13}$ | 0.1266 | 0.07146 | 0.3642 | 0.08205 | 0.04369 |
| $p_{04}$ | -0.3623 | -0.4156 | 1.508 | 0.5399 | -0.01218 |

Although our work uses the methodology described by Cai *et al.* [50] and He *et al.* [5] for handling time delay within LQR, these two design philosophies were discussed for integer order systems. Also, in the original literatures there were no guidelines to optimally choose the LQR weighting matrices. Here we extend these concepts to the context of FO time delay systems with sluggish and oscillatory dynamics using an incommensurate FO state space framework. In addition, a multi-objective trade-off design is proposed to select LQR weighting matrices and the FOPID integro-differential orders and show how this affects two conflicting closed-loop performance objectives.

Also, there have been exhaustive studies on directly tuning the FOPID controller parameters using single or multi-objective optimization algorithms [54]. In Das *et al.* [40], the advantage of weight optimized LQR-FOPID design has been shown over the direct parameter optimized FOPID design. In particular, due to the introduction of the LQR design in the present method, the optimality of state and control trajectories are automatically enforced which is not guaranteed with direct gain and order selection of FOPID controller with a global optimization algorithm [14]. Although the state *vs.* control optimality hold for a particular choice of $Q$ and $R$ due to LQR, yet it is not sufficient from the perspective of traditional time domain performance measures of the control system like set-point tracking, disturbance rejection etc. The dynamic characteristics of the closed loop system keeps on changing with different choice of the LQR weighting matrices. Also, the analytical state *vs.* control optimality is only valid for a particular choice of the weights which motivates the MOO based indirect approach of tuning the LQR weights along with FOPID orders in order to finally obtain the optimal FOPID gains.

In the present study, the design has been restricted considering no change in set-point ($r = 0$) to develop the optimal state feedback stabilization scheme for the incommensurate FO system using LQR. Then a unit step change in the set-point was applied to select the LQR weights and FOPID orders through MOO and the disturbance rejection performances are also found satisfactory. In general LQR theory cannot automatically enforce a disturbance rejection criteria, although in process control it is considered as an essential criterion. Our design is mainly based on LQR to preserve the state optimality of the FOPID control loop and in addition the tracking and control effort were optimized for a unit set-point change. The disturbance rejection performance was not optimized, since with such a criterion the resulting controller parameters would have been totally different and it may not always produce good set-point tracking and low controller effort [54]. The reported simulation study for both the oscillatory and sluggish FO plants show that with the proposed method apart from the considered objectives – ITSE and ISDCO, also a good disturbance rejection performance is achieved. Earlier literatures like [15], [54] reported FOPID design considering both tracking and disturbance inputs but without the consideration of LQR optimality. Bridging the disturbance rejection criteria with the LQR framework for FOPID design needs further investigation and we leave it as the scope of future research.

## 5. Conclusion

An LQR based improved fractional order $PI^\lambda D^\mu$ controller tuning has been proposed in this paper with optimal selection of weighting matrices for handling FO process with time delay, in a compact NIOPTD template. The optimal choice of the weighting matrices along with the FO differ-integrals of the $PI^\lambda D^\mu$ controller have been obtained through multi-objective NSGA-II algorithm, based on simultaneous minimization of two conflicting time domain integral performance indices – ITSE and ISDCO. Thus, the proposed method preserves the state optimality of LQR and at the same time gives a low error index in the closed loop time response while also ensuring stability and efficiently handling the time delay terms of FO process. These improvements enable the control designer to obtain satisfactory closed loop response while also enjoying the benefits of LQR in the optimal $PI^\lambda D^\mu$ controller tuning. The MOO results in a range of controller parameters lying on the Pareto front as opposed to a single controller obtained by commonly adopting single-objective optimization framework, by satisfying different conflicting time domain objectives. It is shown that there exists a trade-off between the two time domain objectives and an improvement in one performance index would invariably result in a deterioration of the other. Thus the designer can choose a controller according to the specific requirements of his control problem. Our simulation results show that the proposed techniques works well even for a highly oscillatory and a highly sluggish FO system with time delay yielding a range of solutions on the Pareto front. For delay dominant plants our simulation shows He's method and for balanced lag and delay plants Cai's method perform better, whereas for lag-dominant systems the solutions are comparable. Tuning rules for the five optimal LQR-FOPID knobs have been provided as a function of process parameter – delay to lag ratio (*L/T*) and fractional exponent of the process (*α*). Future scope of work may include multi-objective LQR based FOPID controller tuning for unstable and integrating fractional order systems with time delay and extending the concept of FO-LQR to noisy processes using FO Kalman filter and Linear Quadratic Gaussian (LQG) technique.

## Appendix

Median solutions of the non-dominated Pareto fronts for all the processes under investigation have been reported in the supplementary material. This data has been used for the tuning rule generation.

## Reference

[1]  B. Anderson and J. B. Moore, *Optimal control: linear quadratic methods*. Prentice-Hall, Inc., 1990.

[2]   A. O'Dwyer, *Handbook of PI and PID controller tuning rules*, vol. 2. World Scientific, 2009.

[3]   R. Vilanova and A. Visioli, *PID Control in the Third Millennium, Advances in Industrial Control*. Springer, New York, 2012.

[4]   G.-R. Yu and R.-C. Hwang, "Optimal PID speed control of brush less DC motors using LQR approach," in *Systems, Man and Cybernetics, 2004 IEEE International Conference on*, vol. 1, 2004, pp. 473–478.

[5]   J.-B. He, Q.-G. Wang, and T.-H. Lee, "PI/PID controller tuning via LQR approach," *Chemical Engineering Science*, vol. 55, no. 13, pp. 2429–2439, 2000.

[6]   I. Podlubny, "Fractional-order systems and PI$\lambda$D$\mu$-controllers," *Automatic Control, IEEE Transactions on*, vol. 44, no. 1, pp. 208–214, 1999.

[7]   C. A. Monje, Y. Chen, B. M. Vinagre, D. Xue, and V. Feliu, *Fractional-order Systems and Controls: Fundamentals and Applications*. London: Springer-Verlag London, 2010.

[8]   S. Das, *Functional fractional calculus*. Springer, 2011.

[9]   S. Das, S. Saha, S. Das, and A. Gupta, "On the selection of tuning methodology of FOPID controllers for the control of higher order processes," *ISA Transactions*, vol. 50, no. 3, pp. 376–388, 2011.

[10]  D. Valerio and J. S. Da Costa, *An introduction to fractional control*, vol. 91. IET, 2013.

[11]  S. Das, I. Pan, S. Das, and A. Gupta, "Improved model reduction and tuning of fractional-order PI$\lambda$D$\mu$ controllers for analytical rule extraction with genetic programming," *ISA Transactions*, vol. 51, no. 2, pp. 237–261, 2012.

[12]  Y. Luo and Y. Chen, *Fractional Order Motion Controls*. John Wiley & Sons, 2012.

[13]  I. Pan and S. Das, "Frequency domain design of fractional order PID controller for AVR system using chaotic multi-objective optimization," *International Journal of Electrical Power & Energy Systems*, vol. 51, pp. 106–118, 2013.

[14]  I. Pan and S. Das, *Intelligent fractional order systems and control*. Springer, 2013.

[15]  S. Das, I. Pan, and S. Das, "Performance comparison of optimal fractional order hybrid fuzzy PID controllers for handling oscillatory fractional order processes with dead time," *ISA Transactions*, vol. 52, no. 4, pp. 550–566, 2013.

[16]  S. Das, A. Gupta, and S. Das, "Generalized frequency domain robust tuning of a family of fractional order PI/PID controllers to handle higher order process dynamics," *Advanced Materials Research*, vol. 403, pp. 4859–4866, 2012.

[17]  S. Das, N.-D. Molla, I. Pan, A. Pakhira, and A. Gupta, "Online identification of fractional order models with time delay: An experimental study," in *Communication and Industrial Application (ICCIA), 2011 International Conference on*, 2011, pp. 1–5.

[18]  S. Das, S. Das, and A. Gupta, "Fractional order modeling of a PHWR under step-back condition and control of its global power with a robust PI$\lambda$D$\mu$ controller," *Nuclear Science, IEEE Transactions on*, vol. 58, no. 5, pp. 2431–2441, 2011.

[19] S. Das, I. Pan, K. Halder, S. Das, and A. Gupta, "LQR based improved discrete PID controller design via optimum selection of weighting matrices using fractional order integral performance index," *Applied Mathematical Modelling*, vol. 37, no. 6, pp. 4253–4268, 2013.

[20] T. Wang, Q. Wang, Y. Hou, and C. Dong, "Suboptimal controller design for flexible launch vehicle based on genetic algorithm: selection of the weighting matrices Q and R," in *Intelligent Computing and Intelligent Systems, 2009. ICIS 2009. IEEE International Conference on*, vol. 2, 2009, pp. 720–724.

[21] J. Liu and Y. Wang, "Design approach of weighting matrices for LQR based on multi-objective evolution algorithm," in *Information and Automation, 2008. ICIA 2008. International Conference on*, 2008, pp. 1188–1192.

[22] M. B. Poodeh, S. Eshtehardiha, A. Kiyoumarsi, and M. Ataei, "Optimizing LQR and pole placement to control buck converter by genetic algorithm," in *Control, Automation and Systems, 2007. ICCAS'07. International Conference on*, 2007, pp. 2195–2200.

[23] O. P. Agrawal, "A general formulation and solution scheme for fractional optimal control problems," *Nonlinear Dynamics*, vol. 38, no. 1–4, pp. 323–337, 2004.

[24] A. Shafieezadeh, Y. Chen, and K. Ryan, "Fractional order filter enhanced LQR for seismic protection of civil structures," *Journal of Computational and Nonlinear Dynamics*, vol. 3, no. 2, p. 021404, 2008.

[25] C. Tricaud and Y. Chen, "An approximate method for numerically solving fractional order optimal control problems of general form," *Computers & Mathematics with Applications*, vol. 59, no. 5, pp. 1644–1655, 2010.

[26] O. P. Agrawal, "A quadratic numerical scheme for fractional optimal control problems," *Journal of Dynamic Systems, Measurement, and Control*, vol. 130, no. 1, p. 011010, 2008.

[27] R. K. Biswas and S. Sen, "Fractional optimal control problems with specified final time," *Journal of Computational and Nonlinear Dynamics*, vol. 6, no. 2, p. 021009, 2011.

[28] R. K. Biswas and S. Sen, "Free final time fractional optimal control problems," *Journal of the Franklin Institute*, vol. 351, no. 2, pp. 941–951, 2014.

[29] X. W. Tangpong and O. P. Agrawal, "Fractional optimal control of continuum systems," *Journal of Vibration and Acoustics*, vol. 131, no. 2, p. 021012, 2009.

[30] R. K. Biswas and S. Sen, "Fractional optimal control problems: a pseudo-state-space approach," *Journal of Vibration and Control*, vol. 17, no. 7, pp. 1034–1041, 2011.

[31] Y. Ding, Z. Wang, and H. Ye, "Optimal control of a fractional-order HIV-immune system with memory," *Control Systems Technology, IEEE Transactions on*, vol. 20, no. 3, pp. 763–769, 2012.

[32] P. M. Czyronis, "Dynamic Programming Problem for Fractional Discrete-Time Dynamic Systems. Quadratic Index of Performance Case," *Circuits, Systems, and Signal Processing*, pp. 1–19, 2014.

[33] A. Dzielinski and P. Czyronis, "Fixed final time and free final state optimal control problem for fractional dynamic systems-linear quadratic discrete-time case," *Bulletin of the Polish Academy of Sciences: Technical Sciences*, vol. 61, no. 3, pp. 681–690, 2013.

[34] N. Özdemir, O. P. Agrawal, B. B. Iskender, and D. Karadeniz, "Fractional optimal control of a 2-dimensional distributed system using eigenfunctions," *Nonlinear Dynamics*, vol. 55, no. 3, pp. 251–260, 2009.

[35] M. M. Hasan, X. W. Tangpong, and O. P. Agrawal, "Fractional optimal control of distributed systems in spherical and cylindrical coordinates," *Journal of Vibration and Control*, vol. 18, no. 10, pp. 1506–1525, 2012.

[36] T. L. Guo, "The Necessary Conditions of Fractional Optimal Control in the Sense of Caputo," *Journal of Optimization Theory and Applications*, vol. 156, no. 1, pp. 115–126, 2013.

[37] Y. Li and Y. Chen, "Fractional order linear quadratic regulator," in *Mechtronic and Embedded Systems and Applications, 2008. MESA 2008. IEEE/ASME International Conference on*, 2008, pp. 363–368.

[38] D. Sierociuk and B. M. Vinagre, "Infinite horizon state-feedback LQR controller for fractional systems," in *Decision and Control (CDC), 2010 49th IEEE Conference on*, 2010, pp. 6674–6679.

[39] S. Saha, S. Das, S. Das, and A. Gupta, "A conformal mapping based fractional order approach for sub-optimal tuning of PID controllers with guaranteed dominant pole placement," *Communications in Nonlinear Science and Numerical Simulation*, vol. 17, no. 9, pp. 3628–3642, 2012.

[40] S. Das, I. Pan, K. Halder, S. Das, and A. Gupta, "Optimum weight selection based LQR formulation for the design of fractional order PIλDμ controllers to handle a class of fractional order systems," in *Computer Communication and Informatics (ICCCI), 2013 International Conference on*, 2013, pp. 1–6.

[41] R. Caponetto, G. Dongola, L. Fortuna, and A. Gallo, "New results on the synthesis of FO-PID controllers," *Communications in Nonlinear Science and Numerical Simulation*, vol. 15, no. 4, pp. 997–1007, 2010.

[42] S. Das, I. Pan, and K. Sur, "Artificial neural network based prediction of optimal pseudo-damping and meta-damping in oscillatory fractional order dynamical systems," in *Advances in Engineering, Science and Management (ICAESM), 2012 International Conference on*, 2012, pp. 350–356.

[43] K. Deb, *Multi-objective optimization using evolutionary algorithms*, vol. 16. John Wiley & Sons, 2001.

[44] N. Aguila-Camacho, M. A. Duarte-Mermoud, and J. A. Gallegos, "Lyapunov functions for fractional order systems," *Communications in Nonlinear Science and Numerical Simulation*, vol. 19, no. 1, pp. 2951–2957, 2014.

[45] J.-C. Trigeassou, N. Maamri, J. Sabatier, and A. Oustaloup, "A Lyapunov approach to the stability of fractional differential equations," *Signal Processing*, vol. 91, no. 3, pp. 437–445, 2011.

[46] J.-C. Trigeassou, N. Maamri, and A. Oustaloup, "Lyapunov Stability of Linear Fractional Systems: Part 1—Definition of Fractional Energy," in *ASME 2013 International Design Engineering Technical Conferences and Computers and Information in Engineering Conference*, 2013, pp. V004T08A025–V004T08A025.

[47] J. Sabatier, M. Moze, and C. Farges, "LMI stability conditions for fractional order systems," *Computers & Mathematics with Applications*, vol. 59, no. 5, pp. 1594–1609, 2010.

[48] M. Moze, J. Sabatier, and A. Oustaloup, "LMI tools for stability analysis of fractional systems," in *ASME 2005 International Design Engineering Technical Conferences and Computers and Information in Engineering Conference*, 2005, pp. 1611–1619.

[49] J. E. Normey-Rico and E. F. Camacho, *Control of dead-time processes*. Springer, 2007.

[50] G.-P. Cai, J.-Z. Huang, and S. X. Yang, "An optimal control method for linear systems with time delay," *Computers & Structures*, vol. 81, no. 15, pp. 1539–1546, 2003.

[51] B. Basu and S. Nagarajaiah, "A wavelet-based time-varying adaptive LQR algorithm for structural control," *Engineering Structures*, vol. 30, no. 9, pp. 2470–2477, 2008.

[52] A. Ruszewski, "Stability regions of closed loop system with time delay inertial plant of fractional order and fractional order PI controller," *Bulletin of the Polish Academy of Sciences, Technical Sciences*, vol. 56, no. 4, pp. 329–332, 2008.

[53] A. Tepljakov, E. Petlenkov, and J. Belikov, "FOMCON: a MATLAB toolbox for fractional-order system identification and control," *International Journal of Microelectronics and Computer Science*, vol. 2, no. 2, pp. 51–62, 2011.

[54] S. Das and I. Pan, "On the Mixed H2/H∞ Loop Shaping Trade-offs in Fractional Order Control of the AVR System," *Industrial Informatics, IEEE Transactions on*, vol. 10, no. 4, pp. 1982 – 1991, 2014.

# Supplementary Material

Table 1: Median solutions on the non-dominated Pareto front

| Process Type | $K_p$ | $K_i$ | $K_d$ | $\lambda$ | $\mu$ | $L/T$ | $\alpha$ |
|---|---|---|---|---|---|---|---|
| Lag dominant | 0.2742 | 0.2952 | 1.4643 | 0.5908 | 0.0287 | 0.25 | 0.2 |
| | 0.1368 | 0.2376 | 1.1798 | 0.6866 | 0.0138 | 0.25 | 0.4 |
| | 0.1630 | 0.2653 | 1.0717 | 0.8634 | 0.0284 | 0.25 | 0.6 |
| | 0.3392 | 0.3684 | 1.4124 | 0.9979 | 0.0206 | 0.25 | 0.8 |
| | 0.2218 | 0.3013 | 1.7505 | 1.3476 | 0.0152 | 0.25 | 1.2 |
| | 0.7017 | 0.5316 | 1.9313 | 1.2673 | 0.2538 | 0.25 | 1.4 |
| | 0.7197 | 0.5679 | 1.8582 | 1.2067 | 0.5152 | 0.25 | 1.6 |
| | 0.6429 | 0.5225 | 1.8176 | 1.2340 | 0.6094 | 0.25 | 1.8 |
| Balanced lag and delay | 0.2271 | 0.6640 | 0.0095 | 0.4276 | 0.9472 | 1 | 0.2 |
| | 0.4045 | 0.6957 | 0.0333 | 0.4963 | 0.2792 | 1 | 0.4 |
| | 2.3864 | 0.3661 | 3.1882 | 0.7451 | 0.0311 | 1 | 0.6 |
| | 0.4247 | 0.3243 | 1.2120 | 0.9049 | 0.1082 | 1 | 0.8 |
| | 0.8149 | 0.3449 | 1.6216 | 0.9846 | 0.1507 | 1 | 1.2 |
| | 0.7910 | 0.3327 | 1.5229 | 1.0294 | 0.1728 | 1 | 1.4 |
| | 0.7075 | 0.3045 | 1.3405 | 1.0799 | 0.2139 | 1 | 1.6 |
| | 0.7338 | 0.2744 | 1.2289 | 1.1583 | 0.2126 | 1 | 1.8 |
| Delay dominant | 0.3317 | 0.5251 | 0.0461 | 0.6242 | 0.4754 | 4 | 0.2 |
| | 0.3264 | 0.5035 | 0.0558 | 0.7560 | 0.6224 | 4 | 0.4 |
| | 0.3046 | 0.4583 | 0.0496 | 0.8389 | 0.6297 | 4 | 0.6 |
| | 0.2921 | 0.4352 | 0.0494 | 0.8956 | 0.4760 | 4 | 0.8 |
| | 0.2625 | 0.3784 | 0.0442 | 0.9515 | 0.6183 | 4 | 1.2 |
| | 0.2528 | 0.3572 | 0.0416 | 0.9759 | 0.6438 | 4 | 1.4 |
| | 0.2398 | 0.3373 | 0.0402 | 0.9949 | 0.6949 | 4 | 1.6 |
| | 0.2111 | 0.3280 | 0.0406 | 0.9985 | 0.7702 | 4 | 1.8 |